\documentclass[11pt]{amsart}
\usepackage{graphicx, color}
\usepackage{amssymb}
\usepackage{amsmath}
\usepackage{latexsym}
\usepackage{setspace}
\makeatletter
\newcommand{\C}{\mathbb C}
\newcommand{\R}{\mathbb R}
\newcommand{\Z}{\mathbb Z}

\newtheorem{thm}{Theorem}[section]
\newtheorem{lem}[thm]{Lemma}
\newtheorem{pro}[thm]{Proposition}
\newtheorem{dfn}[thm]{Definition}

\begin{document}
\title[Inverse pseudo orbit tracing property for robust diffeomorphisms]{Inverse pseudo orbit tracing property for robust diffeomorphisms}
\author{Manseob Lee}
\address{Manseob Lee, Department of Mathematics, Mokwon University, Daejeon 302-729, Korea.}
\email{lmsds@mokwon.ac.kr}


\maketitle

\begin{abstract}
Let $f:M\to M$ be a diffeomorphism of a compact smooth manifold
$M$. Herein, we demonstrate that (i) if
 $f$ has the $C^1$ robustly
inverse pseudo orbit tracing property on the chain recurrent set
$\mathcal{CR}(f)$, then $\mathcal{CR}(f)$ is hyperbolic of $f$ and
(ii) if $f$ has the $C^1$ robustly inverse pseudo orbit tracing
property on a nontrivial transitive set $\Lambda\subset M$, then
$\Lambda$ is hyperbolic for $f$.

\end{abstract}

\onehalfspace

\section{Introduction}
The inverse pseudo orbit tracing property is a dual notion of the
pseudo orbit tracing property that was introduced by Corless and
Pilyugin \cite{CP}. However, the notions are not the same in
general. Kloeden and Ombach \cite{KO} proved that if an expansive
diffeomorphism $f$ has the pseudo orbit tracing property, then it
has the inverse pseudo orbit tracing property with respect to the
continuous method $\mathcal{T}_h$
 (see the definition in section 2). Regarding Lewowicz's results \cite{Lew}, the
Pseudo-Anosov map $f$ of a compact surface $S$ contains the
inverse pseudo orbit tracing property with respect to the class of
the continuous method $\mathcal{T}_h$; however, it is expansive
and not topologically stable. Therefore, it does not has the
pseudo orbit tracing property. To study the hyperbolic structure
(Anosov, structurally stable, Axiom A, $\Omega$-stable,
hyperbolic, etc.), the pseudo orbit tracing theories are highly
useful concepts. In fact, the concepts are close to the hyperbolic
structure. Robinson \cite{Rob} and Sakai \cite{Sa} proved that a
diffeomorphism $f$ of a compact smooth manifold $M$ belongs to the
$C^1$ interior of the set of all diffeomorphisms exhibiting the
pseudo orbit tracing property if and only if it has the hyperbolic
structure. Pilyugin \cite{Pi1} proved that a diffeomorphism $f$ of
a compact smooth manifold $M$ belongs to the $C^1$ interior of the
set of diffeomorphisms exhibiting the inverse pseudo orbit tracing
property with respect to the continuous method $\mathcal{T}_c$
(see the definition in section 2) if and only if it has the
hyperbolic structure.

Lee  shown in \cite{Lee} that if a diffeomorphism $f$ of a compact
smooth manifold $M$ is topologically stable, then it has the
inverse pseudo orbit tracing property with respect to the class of
the continuous method $\mathcal{T}_d $ (see the definition in
section 2). Bowen \cite{Bo} proved that if a diffeomorphism $f$ of
a compact smooth manifold $M$ is hyperbolic, then it has the
pseudo orbit tracing property. Lee \cite{Lee} proved that if a
diffeomorphism $f$ of a compact smooth manifold $M$ is hyperbolic,
then it has the inverse pseudo orbit tracing property with respect
to the class of the continuous method $\mathcal{T}_d.$ Therefore,
we know that if a diffeomorphism $f$ has the hyperbolic structure,
then it has the pseudo orbit tracing and inverse pseudo orbit
tracing properties with respect to the class of the continuous
method $\mathcal{T}_d.$

However, regarding the local dynamical systems with the $C^1$
robust property (see  definition 3.1), the results of two concepts
are different. Lee \cite{L} proved that if a diffeomorphism $f$
has the $C^1$ robustly pseudo orbit tracing property on the
transitive set $\Lambda$, then $\Lambda$ is a hyperbolic set for
$f$. Lee and Lee \cite{LL} proved that if a diffeomorphism $f$ has
the $C^1$ robustly inverse pseudo orbit tracing property with
respect to the class of the continuous method $\mathcal{T}_d$ on
the transitive set $\Lambda$, then $\Lambda$ admits a dominated
splitting for $f.$ However, it is still unclear if a
diffeomorphism $f$ has the inverse pseudo orbit tracing property
with respect to the class of the continuous method $\mathcal{T}_d$
on the transitive set $\Lambda$, thus causing $\Lambda$ to be
hyperbolic set for $f$.
 Therefore, we will prove the problem herein, which is the primary theorem.

The paper is organized as follows. In section 2, we introduce the
pseudo orbit tracing and inverse pseudo orbit tracing properties.
In section 3, we introduce the basic notions and primary theorems.
In section 4, we prove Theorem A. Finally, in section 5, we prove
Theorem B.

\section{Inverse pseudo orbit tracing property}
  Let $M$ be a
compact smooth Riemannian manifold without boundary, and let ${\rm
Diff}(M)$ be the space of $C^1$ diffeomorphisms of $M$ with the
$C^1$ topology. Let $\Lambda\subset M$ be a closed $f$-invariant
set. For any $\delta>0$,
 a sequence of points $\{x_i: i\in\Z\}\subset \Lambda$ is regarded as the {\it $\delta$ pseudo orbit} of $f$ if $d(f(x_i), x_{i+1})<\delta$
 $\forall i\in\Z.$ We say that a diffeomorpshim $f$ has the {\it pseudo orbit tracing
 property} on $\Lambda$ if for any $\epsilon>0$, we can find $\delta>0$
 such that for any $\delta$ pseudo orbit
 $\{x_i:i\in\Z\}\subset\Lambda$, a point $y\in M$ exists such
 that $d(f^i(y), x_i)<\epsilon$ $\forall i\in\Z.$ If $\Lambda=M$, then we say that a diffeomorphism $f$ has the {\it shadowing property}.
   It is known that a diffeomorphism $f$ has the pseudo orbit tracing
   property if and only if $f^n$ has the pseudo orbit tracing property for
   all $n\in\Z\setminus\{0\}$; further, if $f$ has the pseudo orbit tracing property,
   then $f$ has the pseudo orbit tracing property on $\Lambda.$

   Let $M^{\Z}$ be the space of all two-sided
sequences $\{x_i:i\in\Z\}\subset M$ endowed with the product
topology. For any $\delta>0$, we denote by $\Gamma_f(\delta)$ the
set of all $\delta$ pseudo orbit of $f$. A mapping $\xi:M\to
\Gamma_f(\delta)\subset M^{\Z}$ is regarded as {\it
$\delta$-method} for $f$ if $\xi(x)_0=x,$ and $\xi(x)$ is a
$\delta$ pseudo orbit of $f$ through $x,$  where $\xi(x)_0$ means
that the $0$th component of $\xi(x).$ Herein, we set
$\xi(x)=\{\xi(x)_i:i\in\Z\}.$ We say that $\xi$ is a {\it
continuous $\delta$-method} for $f$ if the map $\xi$ is
continuous. We denote by $\mathcal{T}_0(f, \delta)$ the set of all
$\delta$ methods, and by $\mathcal{T}_{c}(f, \delta)$ the set of
all continuous $\delta$ methods. For a homeomorphism $g:M\to M$
with $d_0(f, g)<\delta,$ $g$ induces a continuous $\delta$ method
$\xi(g)$ for $f$ such that
$$\xi(g)(x)=\{g^n(x):n\in\Z\},$$ where $d_0$ is the $C^0$ metric.
For any $\delta>0$, we denote by $\mathcal{T}_h(f, \delta)$ the
set of all continuous $\delta$ methods $\xi(g)$  for $f$ which are
induced by a homeomorphism $g:M\to M$ with $d_0(f, g)<\delta$.
   According to the notions above, we define a strong continuous method that is induced by diffeomorphisms.
    For any $\delta>0$ and  a diffeomorphism $g:M\to M$ with $d_1(f, g)<\delta,$ $g$ induces a continuous $\delta$ method $\xi(g)$  for $f$ such that $$\xi(g)(x)=\{g^n(x):n\in\Z\},$$ where $d_1$ is the $C^1$ metric.
    For any $\delta>0$,  we denote by $\mathcal{T}_d(f, \delta)$ the set of all continuous $\delta$ methods $\xi(g)$ which is induced by a diffeomorphism $g$ for
       $f$ which are induced by a diffeomorphism $g:M\to M$ with $d_1(f,
       g)<\delta$.
  We set $$\mathcal{T}_{a}(f)=\bigcup_{\delta>0}\mathcal{T}_{a}(f, \delta),$$ where $a=0, c, h, d.$ It is clear that $$\mathcal{T}_d(f)\subset\mathcal{T}_h(f)\subset\mathcal{T}_c(f)\subset\mathcal{T}_0(f).$$

We say that a diffeomorphism $f$ has the {\it $\mathcal{T}_a
$-inverse pseudo orbit tracing property} if for any $\epsilon>0$,
there is $\delta>0$ such that for any $\delta$ method
$\xi\in\mathcal{T}_a(f, \delta)$ and any point $x\in M$, a point
$y\in M$ exists such that
$$d(f^n(x), \xi(y)_n)<\epsilon,$$ for all $n\in\Z,$ where  $a=0,
c, h, d.$

   We say that a diffeomorphism $f$ has the {\it inverse pseudo orbit tracing property} with respect to the class of the methods $\mathcal{T}_a $ if it has the $\mathcal{T}_a$ inverse pseudo orbit tracing
   property, where  $a=0, c, h, d.$

Lee and Park \cite{LPark} proved that for a unit circle $S$, a
diffeomorphism $f:S\to S$ has the pseudo orbit tracing property if
and only if $f$ exhibits the inverse pseudo orbit tracing property
with respect to the class of the continuous method
$\mathcal{T}_h$. Sakai \cite{Sa1} proved that a diffeomorphism $f$
of a compact smooth manifold $M$ belongs to the $C^1$ interior of
the set of diffeomorphisms exhibiting the inverse pseudo orbit
tracing property with respect to the class of the continuous
method $\mathcal{T}_h$ then it has a hyperbolic
structure(structurally stable). It was also proved in \cite{Lee}
that if a diffeomorphism $f$ of a compact smooth manifold $M$
belongs to the $C^1$ interior of the set of diffeomorphisms
exhibiting the inverse pseudo orbit tracing property with respect
to the class of the continuous method $\mathcal{T}_d$ then it has
a hyperbolic structure. We denote by $\mathcal{ISP}_a$ the set of
all diffeomorphisms having the inverse (structurally stable)
property with respect to the class of the methods $\mathcal{T}_{a}
(a=0, c, h, d)$. Let $int\mathcal{ISP}_a$ be the  $C^1$ interior
of the set of all diffeomorphisms having the inverse (structurally
stable) property with respect to the class of the methods
$\mathcal{T}_a (a=0, c, h, d).$ According to the results of
Pilyugin \cite{Pi1}, Sakai \cite{Sa1}, and Lee \cite{Lee},
$$int\mathcal{ISP}_c=int\mathcal{ISP}_h =int\mathcal{ISP}_d.$$ By
definition,  we know that
$\mathcal{ISP}_c\subset\mathcal{ISP}_h\subset\mathcal{ISP}_d.$
However,
$\mathcal{ISP}_c\not=\mathcal{ISP}_h\not=\mathcal{ISP}_d,$ in
general. It is noteworthy that $f$ has the inverse pseudo orbit
tracing property with respect to the class of the continuous
method $\mathcal{T}_d$ if and only if $f^n$ has the inverse pseudo
orbit tracing property with respect to the class of the continuous
method $\mathcal{T}_d$, for all $n\in\Z\setminus\{0\}$ (see
\cite{Lee}). It is clear that if $f$ has the inverse pseudo orbit
tracing property with respect to the class of the continuous
method $\mathcal{T}_d$, then $f$ has the inverse pseudo orbit
tracing property on $\Lambda\subset M$ with respect to the class
of the continuous method $\mathcal{T}_d$.

 In this study, we consider the
 inverse pseudo orbit
tracing property with respect to the class of the
 continuous method $\mathcal{T}_d$. Therefore, we use the following expression: {\it a diffeomorphism $f$ has the inverse pseudo orbit
tracing property}.
 This means that a diffeomorphism $f$ has the
 inverse pseudo orbit
tracing property with respect to the class of the
 continuous method $\mathcal{T}_d$.

\section{Basic notions and Theorems}
 In this section, we introduce some notions and primary theorems. Let
 $M$ be as before, and let $f\in{\rm Diff}(M).$
For any $x\in M$,
 $Orb(x)=\{f^n(x):n\in\Z\}$ denotes the orbit of $x.$ A point $p\in
 M$ is called {\it periodic} if $\pi(p)>0$ such that
 $f^{\pi(p)}(p)=p,$ where $\pi(p)$ is the period of $p.$ We denote by
 $P(f)$ the set of all periodic points of $f.$ A point $x\in M$
 is called {\it nonwandering} if in a neighborhood $U$ of $x$, $n>0$ such that $f^n(U)\cap U\not=\emptyset.$ We denote by
 $\Omega(f)$ the set of all nonwandering points of $f.$ It is
 known that $\overline{P(f)}\subset\Omega(f).$  For given $x, y\in M$,
we write $x\rightsquigarrow y$ if for any $\delta>0$, a
$\delta$-pseudo orbit $\{x_i\}_{i=0}^{n}(n>1)$ of $f$ exists such
that $x_0=x$ and $x_n =y.$  We write $x\leftrightsquigarrow y$ if
$x\rightsquigarrow y$ and $y\rightsquigarrow x.$ The set $\{x\in
M: x\leftrightsquigarrow x\}$ is called the {\it chain recurrent
set} of $f$ and is denoted by $\mathcal{CR}(f).$ It is known that
$\Omega(f)\subset \mathcal{CR}(f),$ and $\mathcal{CR}(f)$ is a
closed $f$-invariant set.

 A closed $f$-invariant set
$\Lambda\subset M$ is called {\it hyperbolic} for $f$ if the
tangent bundle $T_{\Lambda}M$ exhibits a $Df$-invariant splitting
$E^s\oplus E^u$ and constants $C>0$ and $0<\lambda<1$ exist such
that
$$\|D_xf^n|_{E_x^s}\|\leq C\lambda^n\;\;{\rm and}\;\;\|D_xf^{-n}|_{E_x^u}\|\leq C\lambda^{n} $$
for all $x\in \Lambda$ and $n\geq 0.$

 We say
that $f$ satisfies {\it Axiom A} if the nonwandering set
$\Omega(f)$ is hyperbolic and it is the closure of $P(f).$

According to Smale \cite{S}, if $f$ satisfies Axiom A, then the
nonwandering set
$\Omega(f)=\Lambda_1\cup\Lambda_2\cup\dots\cup\Lambda_m,$ where
$\Lambda_i$ are compact, disjoint, invariant sets, and each
$\Lambda_i$ contains dense periodic orbits. The sets $\Lambda_1,
\ldots, \Lambda_m$ are called the {\it basic sets}. For a basic
set $\Lambda_i$, we define the following: $$W^s(\Lambda_i)=\{x\in
M: \lim_{n\to\infty} d(f^n(x), \Lambda_i)=0\}, \ \mbox{and}$$
$$W^u(\Lambda_i)=\{x\in
M: \lim_{n\to-\infty} d(f^n(x), \Lambda_i)=0\}.$$

For the basic sets $\Lambda_i (1\leq i\leq n)$, we define
$\Lambda_i>\Lambda_j$ if $$(W^s(\Lambda_i)\setminus\Lambda_i)\cap
W^u(\Lambda_j)\not=\emptyset.$$ We say that $f$ satisfies the {\it
no-cycle condition} if
$\Lambda_{i_0}>\Lambda_{i_1}>\cdots>\Lambda_{i_j}>\Lambda_{i_0}$
cannot occur among the basic sets.


Let $\Lambda\subset M$ be a closed $f$-invariant set. We say that
 $\Lambda$ is {\it locally maximal} if a neighborhood $U$ of
 $\Lambda$ exists such that $\Lambda=\bigcap_{n\in\Z}f^n(U).$
\begin{dfn} Let $f\in{\rm Diff}(M).$ We say that $f$ has the
$C^1$ robustly $\mathcal{P}$ property on $\Lambda$ if a $C^1$
neighborhood $\mathcal{U}(f)$ of $f$ and a neighborhood $U$ of
$\Lambda$ exist such that (i) $\Lambda=\bigcap_{n\in\Z}f^n(U),$
and (ii) for any $g\in\mathcal{U}(f)$, $g$ has the $\mathcal{P}$
property on $\Lambda_g$ where $\Lambda_g$ is the continuation of
$\Lambda.$
\end{dfn}

In the definition, if $\mathcal{P}$ is the pseudo orbit tracing,
then it was defined by Lee, Moriyasu, and Sakai \cite{LMS}. If
$\mathcal{P}$ is the inverse pseudo orbit tracing, then it was
defined by Lee and Lee \cite{LL}.
 Herein, we use the second case where $\mathcal{P}$ is the
 inverse pseudo orbit
tracing.

It is known that if a closed $f$-invariant set $\Lambda\subset M$
is hyperbolic for $f$, then $f$ has the inverse pseudo orbit
tracing property on $\Lambda$. By the stability of hyperbolic
invariant sets for $f$ (\cite[Theorem 7.4]{Ro}), if a closed
$f$-invariant set $\Lambda$ is hyperbolic for $f$, then a $C^1$
neighborhood $\mathcal{U}(f)$ and a neighborhood $U$ of $\Lambda$
exist such that $\Lambda=\bigcap_{n\in\Z}f^n(U)$; further, for any
$g\in\mathcal{U}(f)$, $\Lambda_g=\bigcap_{n\in\Z}g^n(U)$ is
hyperbolic. Therefore, $g$ has the inverse pseudo orbit tracing
property on
$\Lambda_g.$ Hence, we have the following.\\

\noindent{\bf Theorem A}  {\em Let $f\in{\rm Diff}(M),$ and let
$\mathcal{CR}(f)$ be the chain recurrent set of $f.$ If $f$ has
the $C^1$ robustly inverse pseudo orbit tracing property on
$\mathcal{CR}(f)$,
then $\mathcal{CR}(f)$ is hyperbolic.}\\

A closed $f$-invariant set $\Lambda\subset M$ is called {\it
transitive} for $f$ if a point $x\in\Lambda$ exists such that
$\omega(x)=\Lambda,$ where $\omega(x)$ is the omega limit set of
$x.$ In this study, we consider that a transitive set $\Lambda$ is
nontrivial as it is not one orbit.
 We say that a
closed $f$-invariant set $\Lambda\subset M$ admits a {\it
dominated splitting} for $f$ if the tangent bundle $T_{\Lambda}M$
exhibits a continuous $Df$ invariant splitting $E \oplus F$ and
$C>0$, $0 < \lambda < 1$ such that for all $x \in \Lambda$ and $n
\geq 0$, we have
$$
||Df^n|_{E(x)}|| \cdot ||Df^{-n}|_{F(f^n(x))} || \leq C
\lambda^n.$$ As mentioned in the previous section, if a
diffeomorphism $f$ has the inverse pseudo orbit tracing property
on a transitive set $\Lambda$, then it admits a dominated
splitting for $f$ (see \cite{LL}). According to the results, we
prove the
following.\\

\noindent{\bf Theorem B} {\em Let $f\in{\rm Diff}(M)$ and let
$\Lambda$ be a transitive set of $f.$ If $f$ has the $C^1$
robustly inverse pseudo orbit tracing property
on $\Lambda$, then $\Lambda$ is hyperbolic for $f.$}\\

\section{Proof of Theorem A}

In this section, we prove the hyperbolicity of the chain recurrent
set $\mathcal{CR}(f)$ with the $C^1$ robustly inverse pseudo orbit
tracing property. To prove this, we use a $C^1$ perturbation
lemma, called Franks' lemma. The following is Franks' lemma (see
\cite{F}):

\begin{lem}\label{frank} Let $\mathcal{U}(f)$
be any given $C^1$ neighborhood of $f$. Therefore, $\epsilon>0$
and a $C^1$ neighborhood $\mathcal{U}_0(f)\subset\mathcal{U}(f)$
of $f$ exists such that for a given $g\in \mathcal{U}_0(f)$, a
finite set $S=\{x_1, x_2, \cdots, x_N\}$, a neighborhood $U$ of
$S$, and linear maps $A_i : T_{x_i}M\rightarrow T_{g(x_i)}M$
satisfying $\|A_i-D_{x_i}g\|\leq\epsilon$  for all $1\leq i\leq
N$, there exists $h\in \mathcal{U}(f)$ such that $h(x)=g(x)$ if
$x\in S\cup(M\setminus U)$ and $D_{x_i}h=A_i$ for all $1\leq i\leq
N$.
\end{lem}

Using lemma \ref{frank} and the $C^1$ robustly inverse pseudo
orbit tracing property, an important lemma exists as follows. From
the lemma, we can demonstrate that if a diffeomorphism $f$
exhibits the $C^1$ robustly inverse pseudo orbit tracing property
on $\mathcal{CR}(f)$, then $\mathcal{CR}(f)$ is hyperbolic.

\begin{lem}\label{perhy} Let $\Lambda\subset M$ be a closed $f$-invariant set.
If $f$ has the $C^1$ robustly inverse pseudo orbit tracing
property on $\Lambda$, then for any $g$ $C^1$ close to $f$, every
$p\in\Lambda_g\cap P(g)$ is hyperbolic, where $P(g)$ is the set of
periodic points for $g.$
\end{lem}
\noindent{\bf Proof.} Let $\mathcal{U}(f)$ be a $C^1$ neighborhood
of $f$ and $U$ be a locally maximal neighborhood of $\Lambda$.
Suppose that $g\in\mathcal{U}(f)$ exists such that $g$ contains a
nonhyperbolic periodic point
$p\in\Lambda_g=\bigcap_{n\in\Z}g^n(U).$ Because $p\in\Lambda_g\cap
P(g)$ is not hyperbolic, an eigenvalue $\lambda$ of
$D_pg^{\pi(p)}$ exists such that $|\lambda|=1,$ where $\pi(p)$ is
the period of $p.$
 For simplicity, we may assume that $g^{\pi(p)}(p)=g(p)=p.$ Because $p\in \Lambda_g\cap P(g)$ is not hyperbolic,
an eigenvalue $\lambda$ of $D_pg$ exists such that $|\lambda|=1$.
Therefore, $T_pM=E^c_p\oplus E^s_p\oplus E^u_p$ is the
$D_pg$-invariant splitting of $T_pM$, where $E^{c}_p$ corresponds
to eigenvalues $|\lambda|=1$ of $D_pg$, $E_p^s$ corresponds to
eigenvalues $|\lambda|<1$ of $D_pg$, and $E_p^u$ corresponds to
eigenvalues $|\lambda|>1$ of $D_pg$. According to lemma
\ref{frank}, $g_0$ $C^1$ close to $g$ exists such that
$g_0(p)=g(p)=p$ and $p$ is not hyperbolic for $g_0$. Therefore, we
have only one eigenvalue $\lambda$ of $D_pg_0$ such that
$|\lambda|=1$ and $T_pM=\widetilde{E}_p^c\oplus
\widetilde{E}_p^s\oplus \widetilde{E}_p^u$. If $\lambda\in\R$,
then ${\rm dim}\widetilde{E}_p^c =1$; if $\lambda\in\C$, then
${\rm
dim}\widetilde{E}_p^c =2$.\\

 \noindent{\it Case 1.} Consider
$\lambda\in\R.$ We may assume that
$\lambda=1$(the other case is similar).\\

\noindent Using lemma \ref{frank} again, we obtain $\alpha>0$ with
$B(p, \alpha)\subset U$ and $g_1$ $C^1$ close to $g_0$
$(g_1\in\mathcal{U}(f))$, satisfying
\begin{itemize}
\item[(a)]$g_1(p)=g_0(p)=p$,
\item[(b)]$g_1(x)={\rm exp}_{p}\circ D_{p}g_0\circ
{\rm exp}_{p}^{-1}(x) $ if $x\in B(p, \alpha),$   and
\item[(c)] $g_1(x)=g_0(x)$, if $x\in B(p, 4\alpha).$
\end{itemize}

 We use a nonzero vector
$u\in \widetilde{E}_p^c\subset T_pM$ such that $\|u\|=\alpha/4.$
Subsequently,
$$g_1({\rm exp}_p(u))={\rm exp}_p(D_pg({\rm exp}^{-1}_p({\rm
exp}_p(u)))={\rm exp}_p(u).$$

We set $$\mathcal{J}_p={\rm exp}_p\{t\cdot u:-\frac{\alpha}{4}\leq
t\leq \frac{\alpha}{4}\}.$$ For the small arc $\mathcal{J}_p$, the
following properties hold:
\begin{itemize}
\item[(a)] $\mathcal{J}_p\subset B(p, \alpha)\cap {\rm
exp}_p(\widetilde{E}_p^c(\alpha))$ with the center at $p,$
\item[(b)] $\mathcal{J}_p\subset\Lambda_{g_1}$, and
\item[(c)]
$g_1|_{\mathcal{J}_p}:\mathcal{J}_p\to\mathcal{J}_p$ is the
identity map, \end{itemize} where $\widetilde{E}_p^c(\alpha)$ is
the $\alpha$-ball in $\widetilde{E}_p^c$ centered at the origin
$O_p.$

We denote $\widetilde{E}_p^c=\{u\in T_pM: u_1\not=0,
u_2=\cdots=u_n=0\}$ in the coordinates of the corresponding
neighborhoods. We identify $p$ with $O_p$ and $T_pM$ with $\R^n$
in the coordinates of the corresponding neighborhoods.
Subsequently, we know $p=(0, \cdots,0)$ and
$\widetilde{E}_p^c=\{x\in \R^n: x_1\not=0, x_2=\cdots=x_n=0\}.$
Because $f$ has the $C^1$ robustly inverse pseudo orbit tracing
property on $\Lambda$, $g_1$ has the inverse pseudo orbit tracing
property on $\Lambda_{g_1}=\bigcap_{n\in\Z}g_1^n(U).$ We use
$0<\epsilon<\alpha/16$ and let $0<\delta<\epsilon$ be the number
of inverse pseudo orbit tracing properties for $g_1.$  Given
$\alpha>0$, we define the map
$$g_1|_{B(p, \alpha)}:B(p, \alpha) \to B(p, \alpha)$$ by
$g_1(x)=(x_1, Cx'),$ where $C$ is the hyperbolic part of $D_pg_1$
and $x'=(x_2, x_3, \ldots, x_n).$ We define a diffeomorphism
$h:M\to M$ having the following property,
$$h(x)=\big(x_1+\frac{\delta}{4}, Cx' \big) \ \mbox{and}\ \
h^{-1}(x)=\big(x_1-\frac{\delta}{4}, C^{-1}x' \big),$$ for all
$x=(x_1, x_2, x_3, \ldots, x_n)=(x_1, x') \in B(p, \alpha).$
Therefore, we can obtain a class of the continuous $\delta$ method
$\varphi_h\in\mathcal{T}_d(g_1)$ that is induced by $h$ such that
 $\varphi_h(x)_n=\{h^n(x):n\in\Z\},$ for any $x\in M$. Because $\mathcal{J}_p\subset\Lambda_{g_1}$ and  $g_1$
has the inverse pseudo orbit tracing property on $\Lambda_{g_1}$,
$g_1$ must have the inverse pseudo orbit tracing property on
$\mathcal{J}_p$.

We prove that if
$g_1|_{\mathcal{J}_p}:\mathcal{J}_p\to\mathcal{J}_p$ is the
identity map, then $g_1$ does not have the inverse pseudo orbit
tracing property on $\mathcal{J}_p.$

Firstly,  the pseudo point is in $\mathcal{J}_p.$ Then we have two
cases:  (i). If a pseudo point $y\in \mathcal{J}_p$, then because
$g_1|_{\mathcal{J}_p}:\mathcal{J}_p\to\mathcal{J}_p$ is the
identity map, we can easily demonstrate that $g_1$ does not have
the inverse pseudo orbit tracing property on $\mathcal{J}_p.$
Indeed, we choose $x_0=(2\epsilon, 0, \ldots, 0)\in \mathcal{J}_p$
 such that $d(x_0, p)=2\epsilon$. Because $g_1$ has the inverse
 pseudo orbit
tracing property on $\mathcal{J}_p$, we can use a pseudo point
 $y\in\mathcal{J}_p$ such that $y=p=(0, 0, \cdots, 0)$. Then, we can see that for $n\geq0$,\begin{equation}d(g_1^n(x_0),
 \varphi_h(y)_n)=d(x_0, h^n(y))=d(2\epsilon,
 \frac{\delta}{4})>\epsilon.\end{equation}
Since $g_1$ has the inverse pseudo orbit tracing property on
$\mathcal{J}_p$,  this is a contradiction by (1).

(ii). If a pseudo point $y=(y_1+\delta/4, 0, \cdots,
0)\in\mathcal{J}_p$ with $d(x_0, y)<\epsilon$, then $d(2\epsilon,
y_1+\delta/4)<\epsilon.$ By our construction map $h:M\to M$, $j>0$
exists such that $y_1+(\delta/4)j>3\epsilon.$
 Thus, $j>0$ exists such
that \begin{equation} d(g_1^j(x_0), \varphi_h(y)_j)=d(x_0,
h^j(y))=d(2\epsilon,
y_1+\frac{\delta}{4}j)>\epsilon.\end{equation} According to (2),
$g_1$ does not have the inverse pseudo orbit tracing property on
$\mathcal{J}_p.$

 Therefore, for the chosen point
$x_0\in\mathcal{J}_p$, if a pseudo point $y\in\mathcal{J}_p$, then
$g$ does not have the inverse pseudo orbit tracing property on
$\mathcal{J}_p$.

Finally,  we consider that  the pseudo point $y\in M$ has to
remain in $B(x_0, \epsilon)\setminus \mathcal{J}_p.$

Then for any pseudo point $y\in B(x_0,
\epsilon)\setminus\mathcal{J}_p$, because $g_1$ has the inverse
pseudo orbit tracing property on $\mathcal{J}_p$, the following
inequalities hold:
$$d(g_1^n(x_0), \varphi_h(y)_n)=d(g_1^n(x_0), h^n(y))<\epsilon, \
 \forall n\in\Z.$$ Subsequently, by our defined map $h:M\to M$, for
$z=(z_1, z_2, \ldots, z_n)\in M$, we know that for $n\geq0,$
$$h^n(z)=(z_1+\frac{\delta}{4}n, C^nz'),$$ and  $$h^{-n}(z)=(z_1-\frac{\delta}{4}n, C^{-n}z'),$$ where $z'=(z_2, \cdots,
z_n).$ Therefore, we find that $k>0$ such that
$z_1+(\delta/4)k>3\epsilon$. 
 Thus, $k>0$ exists such that \begin{align*}d(g_1^k(x_0), h^k(z))&=d(x_0,
h^k(z))\\&=d((2\epsilon, 0, \ldots, 0),
(z_1+\frac{\delta}{4}k,C^kz'))\\&\geq d(2\epsilon,
z_1+\frac{\delta}{4}k)>\epsilon.\end{align*}

For the point $y\in  B(x_0, \epsilon)$ with $d(x_0, y)<\epsilon,$
by $g_1$ has the inverse pseudo orbit tracing property on
$\mathcal{J}_p$, the following inequality  $d(g_1^n(x_0),
\varphi_h(y)_n)<\epsilon$ holds, for all $n\in\Z$. However, by the
arguments above, $k>0$ such that $y_1+(\delta/4)k>3\epsilon$.
Thus,
$$d(g_1^n(x_0), \varphi_h(y)_k)=d(x_0, h^k(y))=d(2\epsilon,
y_1+\frac{\delta}{4}k)>\epsilon.$$ Because $g_1$ has the inverse
pseudo orbit tracing property on $\mathcal{J}_p,$ this is a
contradiction. Thus, if
$g_1|_{\mathcal{J}_p}:\mathcal{J}_p\to\mathcal{J}_p$ is the
identity map, then $g_1$ does not have the inverse pseudo orbit
tracing property on $\mathcal{J}_p.$\\

 \noindent{\it Case 2.} Consider
$\lambda\in\C.$ To avoid complexity, we assume that
$g^{\pi(p)}(p)=g(p)=p.$ According to lemma \ref{frank}, $\alpha>0$
exists with $B(p, \alpha)\subset U$ and $g_1$ $C^1$ close to $g$
exhibiting the following properties:
\begin{itemize}
\item[(a)] $g_1(x)={\rm exp}_p\circ D_pg\circ{\rm exp}_p^{-1}(x)$,
if $x\in B(p, \alpha),$
\item[(b)] $g_1(x)=g(x)$, if $x\not\in B(p, 4\alpha)$, and
\item[(c)] $g_1(p)=g(p)=p.$
\end{itemize}
By modifying the map $D_pg_1$, $l>0$ exists such that
$D_pg_1^l(v)=v$ for any $v\in E_p^c(\alpha)\cap {\rm
exp}_p^{-1}(B(p, \alpha)).$ Thus, a small arc
$\mathcal{C}_p\subset {\rm exp}_p(E_p^c(\alpha))\cap B(p, \alpha)$
can be obtained such that $g_1^l(\mathcal{C}_p)=\mathcal{C}_p$ and
$g_1^l|_{\mathcal{C}_p}:\mathcal{C}_p\to\mathcal{C}_p$ is the
identity map. Because $g_1$ has the inverse pseudo orbit tracing
property, it is evident that $g_1^i$ has the inverse pseudo orbit
tracing property for $i\in\Z\setminus\{0\}.$ Let $g_1^l=g_2$.
Therefore, $g_2|_{\mathcal{C}_p}:\mathcal{C}_p\to\mathcal{C}_p$ is
the
identity map. Thus, as in the proof of case 1, a contradiction will be shown. \hfill$\square$\\

We say that a diffeomorphism $f$ is a {\it star} if a $C^1$
neighborhood $\mathcal{U}(f)$ of $f$ exists such that for any
$g\in\mathcal{U}(f)$, every periodic point in $P(g)$ is
hyperbolic. We denote by $\mathcal{F}(M)$ the set of all star
diffeomorphisms. Aoki \cite{A} and Hayashi \cite{H} proved that if
a diffeomorphism $f$ is a star, then $f$ satisfies Axiom A and
no-cycle condition. It is well known that if $f$ satisfies Axiom
A, then $\overline{P(f)}=\Omega(f)=\mathcal{CR}(f)$ (see
\cite{Sh})and the chain recurrent set $\mathcal{CR}(f)$ is upper
semi-continuous, that is, for any neighborhood $U$ of
$\mathcal{CR}(f)$, $\delta> 0$ such that if $d_{C^0}(f, g) <
\delta (g \in{\rm Diff}(M))$, then $\mathcal{CR}(g)\subset U$,
where $d_{C^0}$ is the  $C^0$-metric on
${\rm Diff}(M)$ (see \cite[Corollary 3 (a)]{Hu}).\\

\noindent{\bf Proof of Theorem A.} The arguments above are
sufficient to demonstrate that $f$ is a star. Let $\mathcal{U}(f)$
be a $C^1$ neighborhood of $f$ and a neighborhood $U$ of
$\mathcal{CR}(f).$ Because the chain recurrent set
$\mathcal{CR}(f)$ is upper semi-continuous, we know that
$\mathcal{CR}(g)\subset U$; therefore,
$P(g)\subset\mathcal{CR}(g)\subset\Lambda_g=\bigcap_{n\in\Z}g^n(U).$
 Because $f$ has the $C^1$ robustly inverse
pseudo orbit tracing property on $\mathcal{CR}(f)$, according to
lemma \ref{perhy}, every $p\in\Lambda_g\cap P(g)= P(g)$ is
hyperbolic for any $g\in\mathcal{U}(f)$. Therefore, $f$ is a star,
that is, $f$ satisfies Axiom A and the no-cycle condition. Thus,
the chain recurrent set $\mathcal{CR}(f)$ is hyperbolic.
\hfill$\square$

\section{Proof of Theorem B}
In this section, we introduce a local star condition. Using the
condition, we demonstrate that if a diffeomorphism $f$ exhibits
the $C^1$ robustly inverse pseudo orbit tracing property on a
transitive set $\Lambda$, then $f$ is a star on $\Lambda$.
Therefore, the transitive set $\Lambda$ is hyperbolic for $f.$
Let $\Lambda\subset M$ be a closed $f$-invariant set. We say that
a diffeomorphism $f$ is a {\it star} on $\Lambda$ if a $C^1$
neighborhood $\mathcal{U}(f)$ of $f$ and a neighborhood $U$ of
$\Lambda$ exist such that for any $g\in\mathcal{U}(f)$, every
$p\in \Lambda_g\cap P(g)$ is hyperbolic, where
$\Lambda_g=\Lambda_g(U)=\bigcap_{n\in\Z}g^n(U)$ is the
continuation of $\Lambda.$ It is clear that if $\Lambda=M$, then
$f$ is a star. We denote by $\mathcal{F}(\Lambda)$ the set of all
diffeomorphisms that are stars on $\Lambda$.

\begin{lem}\label{lstar} Let $\Lambda$ be a closed invariant set of $f$. If $f$ exhibits
the $C^1$ robustly inverse pseudo orbit tracing property on
$\Lambda$, then $f\in\mathcal{F}(\Lambda).$
\end{lem}
\noindent{\bf Proof.} Suppose that $f$ exhibits the $C^1$ robustly
inverse pseudo orbit tracing property on $\Lambda$. By the
definition of $\mathcal{F}(\Lambda)$, a $C^1$ neighborhood
$\mathcal{U}(f)$ of $f$ and a neighborhood $U$ of $\Lambda$ exist
such that for any $g\in\mathcal{U}(f)$, every $p\in\Lambda_g\cap
P(g)$ is hyperbolic. Subsequently, the proof is the same as that
of lemma
\ref{perhy}.\hfill$\square$\\

If $p$ is a hyperbolic periodic point, then a $C^1$ neighborhood
$\mathcal{U}(f)$ and a neighborhood $U$ of $p$ exist such that for
any $g\in\mathcal{U}(f)$, a hyperbolic periodic point $p_g\in
P(g)$ exists, where $p_g=\bigcap_{n\in\Z}g^n(U)$ is called {\it
the continuation} of $p.$ Ma\~n\'e \cite[Lemma II.3]{M3} and Lee
and Park \cite[Lemma 2.3]{LP} proved the following:

\begin{pro}\label{mane} Let $\Lambda$ be a transitive set of $f$. Suppose that $f\in\mathcal{F}(\Lambda).$ Therefore, a
$C^1$ neighborhood $\mathcal{U}(f)$ of $f$, constants $C>0,
0<\lambda<1$, and $m \in \mathbb{Z}^+$ exist such that
\begin{itemize}
\item[{\rm(a)}] for each $g \in \mathcal{U}(f)$, if $p$ is a periodic point
of $g$ in $\Lambda_g$ with period $\pi(p, g) \geq m$.
 Therefore,
$$
\prod^{k-1}_{i=0} ||Dg^m|_{E^s(g^{im}(p))}|| < C \lambda^k \quad
\text{and} \quad \prod^{k-1}_{i=0} ||Dg^{-m}|_{E^u(g^{-im}(p))}||
< C \lambda^k,
$$
where $k =[\pi(p,g)/m]$.
\item[{\rm(b)}] $\Lambda$ admits a dominated splitting
$T_{\Lambda}M=E\oplus F$ with ${\rm dim}E={\rm index}(p)$.
\end{itemize}
\end{pro}

 A closed $f$-invariant set $\Lambda\subset M$  is called
an {\it $i$-fundamental limit set} of $f$ if sequences $g_n\to f$
exist as $n\to\infty$ and periodic orbits $P_n$ of $g_n$ with
index $i$ exist such that $\Lambda$ is the Hausdorff limit of
$P_n$. It is noteworthy that the fundamental $i$-limit $\Lambda$
of $f$ is $f$-invariant \cite{Liao}.

\begin{lem}\label{index2}
Let $\Lambda$ be a transitive set and $f\in\mathcal{F}(\Lambda)$.
Then there exist a $C^1$ neighborhood $\mathcal{V}(f)$ of $f$ and
a neighborhood $V$ of $\Lambda$ such that for any integer $i$, if
$g\in\mathcal{U}(f)$ exists such that $g$ exhiits a hyperbolic
periodic point $q\in U$ of index $i$, then $f$ also exhibits a
hyperbolic periodic point of index $i$ in $\Lambda$ and $\Lambda$
is an $i$-fundamental limit set, where $\mathcal{U}(f)$ and $U$
are as the definition of $f\in\mathcal{F}(\Lambda)$.
\end{lem}

\noindent {\bf Proof.} Set $V\subset\overline{V}\subset U$ with an
open neighborhood of $\Lambda$. Let $\mathcal{U}(f)$ be a
neighborhood of $f$ with the following properties: ($a$) for any
$g\in\mathcal{V}(f)(\subset\mathcal{U}(f))$, a continuous path
$\{F_t:0\leq t\leq 1\}\subset {\rm Diff}(M)$ connecting $f$ and
$g$ exists such that any $F_t$ contains no nonhyperbolic periodic
orbits in the neighborhood $V$ of $\Lambda$, ($b$) for any
$g\in\mathcal{V}(f)$, $\bigcap_{i\in
\mathbb{Z}}g^i(\overline{V})=\Lambda_g(U)=\Lambda_g$. We assume
that $g\in\mathcal{U}(f)$ exists such that $g$ contains a
hyperbolic periodic point $q\in U$ of index $i$. Subsequently, we
consider a continuous path $\{F_t:0\leq t\leq1 \}\subset {\rm
Diff}(M)$ connecting $g$ and $f$ such that any $F_t$ contains no
nonhyperbolic periodic orbit in the neighborhood $V$ of $\Lambda$.
If $f$ contains no hyperbolic periodic orbits of index $i$ in
$\Lambda$, then a time $t_0$ exists such that the hyperbolic
periodic orbits of index $i$ is vanished. Without loss of
generality, let $t_0$ be the first time. Therefore, we know that
$F_{t_0}$ contains a nonhyperbolic periodic orbit in $U$; this
contradicts with the path choice. Hence, $f$ also contains a
hyperbolic periodic point of index $i$ in $\Lambda$.

Let $P\subset\Lambda$ be a hyerperbolic periodic orbit of $f$ with
index $i$. By the standard arguments of the connecting lemma (for
instance, see Lemma 2.2 of \cite{GW}), we can apply an arbitrarily
small perturbation $g$ of $f$ such that a  homoclinic orbit $
Orb(x)$ exists with respect to $P$ in $U$, such that the closure
of $ Orb(x)$ is arbitrarily close to the set $\Lambda$ (in
Hausdorff metrics). Applying another arbitrarily small
perturbation if necessary, we can assume that $x$ is a transversal
homoclinic point of $P$. Subsequently, by the shadowing lemma of
hyperbolic set $ Orb(x)\cup P$, we can obtain hyperbolic periodic
orbits of $g$ with index $i$ That are arbitrarily close to $
Orb(x)\cup P$, and hence close to $\Lambda$. This ends the proof
of the second part of lemma \ref{index2}.
\hfill$\square$\\

For any $f\in {\rm Diff}(M)$ and $x\in M$, we denote
$$D^s(x)=D^s(x,f)=\{v\in T_xM: \|Df^n(v)\|\to 0 \ {\rm as}\
n\to+\infty\},$$
$$D^u(x)=D^u(x,f)=\{v\in T_xM: \|Df^n(v)\|\to 0 \ {\rm as}\  n\to-\infty\}.$$

In \cite{Man1}, a characterization of hyperbolicity is detailed as
Follows:
\begin{pro}
A closed $f$-invariant set $\Lambda\subset M$  is hyperbolic if
and only if $T_xM=D^s(x)\oplus D^u(x)$ for any $x\in \Lambda.$
\end{pro}

A point $x\in M$ without the property $T_xM=D^s(x)\oplus D^u(x)$
is called a {\it resisting point}. A compact $f$-invariant set $K$
is called a {\it minimally nonhyperbolic set} if $K$ is
nonhyperbolic and every compact $f$-invariant proper subset of $K$
is hyperbolic. In \cite {Liao}, minimally nonhyperbolic sets are
divided into two types. If a resisting point $a$ exists in a
minimally nonhyperbolic set $K$ such that $\omega(a)$ and
$\alpha(a)$ are all proper subsets of $K$, then $K$ is called the
{\it simple type}. Otherwise, the nonhyperbolic set is called the
{\it nonsimple type}.

\subsection{Non-existence of heterodimensional cycle} In this
Section, we prove the following proposition: no heterodimensional
cycle exists near $\Lambda$ for the system close to $f$.

\begin{pro}\label{nohetero}
Let $\Lambda$ be a transitive set and $f\in\mathcal{F}(\Lambda)$.
Therefore, a $C^1$ neighborhood $\mathcal{U}(f)$ of $f$ and a
neighborhood $U$ of $\Lambda$ exist such that for any
$g\in\mathcal{U}(f)$, $g$ has no a heterodimensional cycles in
$U$.
\end{pro}

\noindent{\bf Proof.} To derive a contradiction, we may assume
that hyperbolic periodic points $p, q$ exist with different
indices and $x\in W^s(p)\cap W^u(q), y\in W^u(p)\cap W^s(q)$ such
that $Orb(p)\cup Orb(q)\cup Orb(x)\cup Orb(y)\subset U.$ We denote
by $K=Orb(p)\cup Orb(q)\cup Orb(x)\cup Orb(y)$ and $k$ the index
of $p$ and $l$ the index of $q$. Without loss of generality, we
can
assume that $p,q$ are fixed points of $f$ and $k<l$.\\

 A point
$x\in M$ is $C^1$ {\it preperiodic} if for any $C^1$ neighborhood
$\mathcal{U}(f)$ of $f$ and any neighborhood $U$ of $x$,
$g\in\mathcal{U}(f)$ and $y\in U$ exist such that $y$ is a
periodic of $g$. We denote by $P_* (f)$ the set of $C^1$
preperiodic points of $f.$ A point $x\in M$ is called an {\it
$i$-preperiodic} of $f$ ($0\leq i\leq{\rm dim}M$) if for any $C^1$
neighbrohood $\mathcal{U}(f)$ of $f$ and any neighborhood $U$ of
$x$, $g\in\mathcal{U}(f)$ and $y\in U$ exist such that $y$ is a
hyperbolic periodic point of $g$ of index $i$ (see \cite{W1}).

\begin{lem}\label{index1}
$K$ is contained in the $k,l$-fundamental limits of $f$.
Precisely, $g_n\to f$ exists with hyperbolic periodic orbits $p_n$
of index $k$, such that $K$ is the Hausdorff limit $p_n$.
Similarly, $g'_n\to f$ exists with hyperbolic periodic orbits
$q_n$ of index $l$, such that $K$ is the Hausdorff limit $q_n$.
\end{lem}

\noindent{\bf Proof.} Because $x\in W^s(p)\cap W^u(q), y\in
W^u(p)\cap W^s(q)$, for any neighborhoods $U_x$ of $x$, $U_y$ of
$y$, and $U_q$ of $q$, one can obtain a point $z$ with integers
$i_1<i_2<i_3$ such that $f^{i_1}(z)\in U_y, f^{i_2}(z)\in U_q$ and
$f^{i_3}(z)\in U_x$ by Palis' $\lambda$-lemma. By small
perturbations, we can create jumps near $x$ and $y$ such that $z$
is a transversal homoclinic point of $p$ for a diffeomorphism $g$
close to $f$. Because the intersection is transversal, we know
that the set $Org(z,g)\cup  Orb(p,g)$ is a hyperbolic set. By the
pseudo orbit tracing lemma, a hyperbolic periodic orbit $p'$ of
$g$ with the same index of $p$ exists such that it is arbitrarily
close to the set $Org(z,g)\cup  Orb(p,g)$. By choosing
sufficiently small $U_x, U_q$, and $U_y$, we can cause the set
$Org(z,g)\cup Orb(p,g)$ to be arbitrarily close to $K$. This
proves that $K$ is the $k$-preperiodic set of $f$. Similarly, we
can prove that $K$ is the $l$-preperiodic set of $f$. This ends
the proof of lemma
\ref{index1}.  \hfill$\square$\\

Let us consider a sequence of periodic pseudo orbits.
\begin{lem}\label{periodic}  Set any small $\delta>0$ and
$x_p\in Orb^+(x), y_p\in Orb^-(y), x_q\in Orb^-(x)$ and $y_q\in
Orb^+(y)$ with $x_p, y_p\in B(\delta, p), x_q, y_q\in B(\delta,
q)$. Subsequently, for any $\epsilon>0$, $L>0$ such that for any
$n\geq L$, $p_n, q_n$ exist with the following properties
\begin{itemize}
\item[(a)] $p_n, f(p_n), \ldots, f^n(p_n)\in B(\delta, p), \ q_n,
f(q_n), \ldots, f^n(q_n)\in B(\delta, q)$,
\item[(b)] $d(x_p, p_n)\leq\epsilon, d(f^n(p_n), y_p)\leq\epsilon,$
$d(f^n(q_n), x_q)\leq\epsilon$ and $d(y_q, q_n)\leq\epsilon.$
\end{itemize}
\end{lem}
\noindent{\bf Proof.} Let $\delta>0$ and $x_p\in Orb^+(x), y_p\in
Orb^-(y), x_q\in Orb^-(x)$ and $y_q\in Orb^+(y)$ with $x_p, y_p\in
B(\delta, p), x_q, y_q\in B(\delta, q)$. By the inclination lemma
of Palis, $$f^n(B(\delta, p))\to W^u(p)\ \mbox{and} \
f^n(B(\delta, q))\to W^u(q)$$ as $n\to\infty.$ Subsequently, for
any $\epsilon>0$, $L>0$ such that for any $n\geq L$, we can use
$p_n $ and $q_n$ such that
\begin{itemize}
\item[(a)] $p_n, f(p_n), \ldots, f^n(p_n)\in B(\delta, p), \ q_n,
f(q_n), \ldots, f^n(q_n)\in B(\delta, q)$,
\item[(b)] $d(x_p, p_n)\leq\epsilon, d(f^n(p_n), y_p)\leq\epsilon,$
$d(f^n(q_n), x_q)\leq\epsilon$ and $d(y_q, q_n)\leq\epsilon.$
\end{itemize}   \hfill$\square$

Consider an $\epsilon$-pseudo orbit $\mathcal{PO}(m, n)$ for $m,
n>L.$ \begin{align*}\mathcal{PO}(m, n)=\{x, f(x), \ldots,
f^{-1}(x_p), p_n, f(p_n), \ldots, f^{n-1}(p_n), y_p, f(y_p),
\ldots,\\
f^{-1}(y_q), q_m, \ldots, f(q_m), \ldots, f^{m-1}(q_m), x_q,
\ldots, f^{-1}(x), x\}.\end{align*}
\begin{lem} Set any small $\delta>0$; $\epsilon>0$
and $N>L$ exist such that if $n\geq N, m\geq N$, then $g$ $C^1$
close to $f$ exists such that $\mathcal{PO}(m, n)$ is a periodic
orbit of $g$.
\end{lem}
\noindent{\bf Proof.} Let any small $\delta>0$ be fixed and $N>L$.
Because $\mathcal{PO}(m, n)$ is a periodic $\epsilon$-pseudo orbit
of $f$, for some $0<\epsilon\leq \delta$, we can create four small
perturbations in a neighborhood of $\{x_p, x_q, y_p, y_q\}$.
Subsequently, the pseudo orbit $\mathcal{PO}(m, n)$ can be a
periodic orbit for the perturbation.\hfill$\square$

\begin{lem}\label{index} If $\delta>0$ is sufficiently small, then for a fixed $n$,
the index of $\mathcal{PO}(m, n)$(with respect to $g$) will equal
to the index of $q$ as $m$ becomes sufficiently large.
\end{lem}
\noindent{\bf Proof.} From Proposition \ref{mane} and  Lemma
\ref{index1}, we know that the set $K$ contains a dominated
splitting $T_K M=E\oplus F$ with ${\rm dim} E=l$. Because $g$ can
be chosen arbitrarily close to $f$ and $\mathcal{PO}(m,n)$
arbitrarily close to $K$, the dominated splitting can continue for
the periodic orbit $\mathcal{PO}(m,n)$ with respect to $g$.
Without loss of generality, we still use $E\oplus F$ to denote the
dominated splitting. Because $x_q$ is close to $q$, we know that
$Dg|_{E(x_q)}$ is close to $Df|_{E^s(q)}$. By the contraction of
$Df|_{E^s(q)}$, after an easy calculation, we find  that
$E|_{\mathcal{PO}(m,n)}$ is contracting with respect to $g$ if $m$
is sufficiently large. Similarly, $F|_{\mathcal{PO}(m,n)}$ is
expanding if $m$ is sufficiently large. This proves that the
periodic orbit $\mathcal{PO}(m, n)$ of $g$ contains an index equal
to $l$. This ends the proof of the lemma. \hfill$\square$

\bigskip
Now, we can complete the proof of Proposition \ref{nohetero}. We
set $m_0,n_0$ to be sufficiently large. By Lemma \ref{index}, we
know that $m>m_0$ exists such that the index of
$\mathcal{PO}(m,n_0)$ is equal to $l$. Subsequently, we set $m$
and increase $n$. In this process, the index of
$\mathcal{PO}(m,n)$ decreases as $n$ increases. If $m_0,n_0$ is
chosen sufficiently large, we can find that $n>n_0$ such that
$\mathcal{PO}(m,n)$ contains the index $k+1$ and
$\mathcal{PO}(m,n+1)$ contains index $k$. By an easy calculation,
we know that if $m_0,n_0$ is sufficiently large, then
$\mathcal{PO}(m,n)$ must contain an eigenvalue $\lambda$ such that
$|\lambda|^{\frac{1}{\pi(\mathcal{PO}(m,n))}}$ is close to $1$.
This is a contradiction because the set $\Lambda$ satisfies the
local star condition. \hfill$\square$

\subsection{Hyperbolicity of local star transitive sets}
In this section, we will prove that if $f$ satisfies the local
star condition, i.e., the transitive $\Lambda$, then it is
hyperbolic. Assume that $\Lambda$ is not a hyperbolic set for $f$.
By Zorn's lemma, we know that a minimally nonhyperbolic set
$K\subset \Lambda$ exists.

\begin{pro} \label{nonsimple} $K$ cannot be a nonsimple-type minimally nonhyperbolic set.
\end{pro}

\noindent{\bf Proof.} Assume that $K$ is a nonsimple-type
minimally nonhyperbolic set. Without loss of generality, we assume
that a resisting point $a$ exists such that $K=\omega(a)$. Let
$k={\rm min}\{i: \text{there is a } i$-$\text{fundamental limit
set contained in } K\}$. From Proposition \ref{mane} and Lemma
\ref{index1}, we know that a dominated splitting $T_KM=E\oplus F$
exists with ${\rm dim}E=k$. Therefore, by ergodic closing lemma
\cite{M3}, we know that $E$ is contracting.

Now, let $$G=\{x\in
K:\limsup_{n\to+\infty}\frac{1}{n}\sum_{i=0}^{n-1}\log(\|Df^m|_{F(f^{im}(x))}\|)\geq
-\log\lambda\}$$ where $m,\lambda$ are the constants in
Proposition \ref{mane}. It is obvious that $\overline{G}$ is a
nonempty compact invariant subset of $K$.

\bigskip

\noindent {\bf Claim. } {\it $\overline {G}=K$.}

\bigskip

\noindent{\bf Proof of Claim.} Assume $\overline{G}$ is a proper
subset of $K$. Subsequently, we know that $\overline{G}$ is
hyperbolic because $K$ is a minimally nonhyperbolic set. It is
easy to verify that $E\oplus F$ restricted on $\overline{G}$ is
only the hyperbolic splitting over $\overline{G}$.

Because $K=\omega(a)$, we know that $a\notin \overline{G}$. One
can apply a small neighborhood $W$ of $\overline{G}$ such that
$a\notin \overline{W}$ and the locally maximal invariant set in
$W$ is hyperbolic. Because $a\notin \overline{W}$ and
$\overline{G}\subset \omega(a)$, we can obtain a point $b\in K$
such that $b\in \overline{W}\setminus f(W)$ and $ Orb^+(b)\subset
\overline {W}$. We know that $b\notin\overline{G}$. Therefore, we
can obtain
$$\limsup_{n\to+\infty}\frac{1}{n}\sum_{i=0}^{n-1}\log(\|Df^m|_{F(f^{im}(b))}\|)<
-\log\lambda.$$ Let $\{n_i\}$ be a sequence of positive integers
such that $f^{n_im}(b)\to c\in\omega(b)$ as $i\to \infty$.
Subsequently, we can apply $1>\lambda'>\lambda$ and $n_s>n_t$ with
$s,t$ arbitrarily large such that
$$\frac{1}{n_s-n_t}\sum_{i=n_s}^{n_t-1}\log(\|Df^m|_{F(f^{im}(b))}\|)<
-\log\lambda'.$$ Subsequently, by the pseudo orbit tracing
property of the hyperbolic sets, we can obtain a hyperbolic
periodic point $p$ with an arbitrarily large period that traces
the orbit segment
$$\{f^{n_sm}(b), f^{(n_s+1)m}(b),\cdots, f^{(n_t-1)m}(b), f^{n_tm}(b)\}$$ such that
$$\frac{1}{\pi(p)}\sum_{i=0}^{\pi(p)-1}\log(\|Df^m|_{E^u(f^{im}(p))}\|)<
-\log\lambda.$$ This contradicts with Proposition \ref{mane}. This
ends the proof of claim.\\

Further, $K$ is shown as a hyperbolic set by the following
conclusion proven in \cite{Man2}. This contradicts that $K$ is a
nonhyperbolic set of $f$. This ends the proof of Proposition
\ref{nonsimple}. \hfill$\square$

\begin{thm}\cite{Man2}
Let $K$ be a compact invariant set of $f$ and assume that $f$ is a
local star in the neighborhood $U$ of $K$. If a dominated
splitting $T_KM=E\oplus F$ exists with the following two
properties:
\begin{itemize}
\item[(a)] $E$ is contracting, and \item[(b)] constants
$m\in\mathbb{N}$ and $\lambda\in(0,1)$, and a dense subset
$G\subset\Lambda$ exist such that for any $x\in G$,
$$\limsup_{n\to+\infty}\frac{1}{n}\sum_{i=0}^{n-1}\log(\|Df^m|_{F(f^{im}(x))}\|)\geq
-\log\lambda,$$
\end{itemize}
then $F$ is expanding and $K$ is hyperbolic.
\end{thm}
\begin{pro}\label{hetero1} If $K$ is a simple-type minimally nonhyperbolic set of $f$, then $g$ $C^1$ close to $f$ exists such that $g$ has a heterodimensional
cycle in $U$.
\end{pro}

\noindent{\bf Proof.} Let $a$ be a resisting point such that
$\omega(a)$ and $\alpha(a)$ are both the proper subsets of $K$.
From the definition of a minimally nonhyperbolic set, we know that
$K=\omega(a)\cup  Orb(a)\cup \alpha(a)$ and both $\omega(a)$ and
$\alpha(a)$ are hyperbolic sets.

\vskip 0.3cm

\noindent {\bf Claim. } {\it The index of $\omega(a)$ and
$\alpha(a)$ are different.}

\vskip 0.3cm

\noindent {\bf Proof of Claim.} Assume that the index of
$\omega(a)$ and $\alpha(a)$ are same. We denote by $i$ the index
of $\omega(a)$. Subsequently, by the pseudo orbit tracing lemma of
the hyperbolic sets, we know that $\Lambda$ contains hyperbolic
periodic points with index $i$. From Lemma \ref{index2}, we know
that $\Lambda$ is an $i$-fundamental limit. From Proposition
\ref{mane}, we know that $\Lambda$ contains a dominated splitting
$T_\Lambda M=E\oplus F$ with ${\rm dim} E=i$. One can easily
verify that $E(x)=D^s(x)$ and $F(x)=D^u(x)$. This contradicts with
$x$ being a resisting points. This ends the proof of claim.

\vskip 0.3cm

We denote by $i$ the index of $\omega(a)$ and $j$ the index of
$\alpha(a)$. Let $W_1$ be a small neighborhood of $\omega(a)$ such
that the maximal invariant set in $W_1$ is hyperbolic and any two
periodic orbits in $W_1$ are homoclinically related. Let $W_2$ be
a small neighborhood of $\alpha(a)$ such that the maximal
invariant set in $W_1$ is hyperbolic and any two periodic orbits
in $W_2$ are homoclinically related. We can small $W_1,W_2$ such
that $W_1\cap W_2=\emptyset$ and $\Lambda\setminus(W_1\cup
W_2)\neq\emptyset$. Let $P$ be a hyperbolic periodic orbit in
$W_1$, and $Q$ be a hyperbolic periodic orbit in $W_2$. By the
standard argument of connecting lemma, we can perform a
perturbation $g$ such that $g=f$ in $W_1\cup W_2\cup  Orb(a)$ and
$W^u(P,g)\cap W^s(Q,g)\neq\emptyset$. It is noteworthy that $g=f$
in $W_1\cup W_2\cup Orb(a)$, we also have
$\omega(a,g)=\omega(a,g)$ and $\alpha(a,g)=\alpha(a,f).$

\begin{lem}\label{diff}
$a\in \overline{W^s(P,g)}\cap\overline{W^u(Q,g)}$.
\end{lem}

\noindent{\bf Proof.}  For an arbitrarily small $\delta>0$, we can
apply $b\in\omega(a)$ and $n\in \mathbb{Z}$ such that
$d(f^n(a),b)<\delta$; subsequently, we can construct a
$\delta$-pseudo orbit as $$\{\cdots, f^{-2}(b),f^{-1}(b), f^n(a),
f^{n+1}(a), \cdots\}.$$ By the pseudo orbit tracing property of
the hyperbolic set, we can find $y\in W_1$ such that the orbit of
$y$ traces the pseudo orbit. If $\delta$ is sufficiently small, we
can obtain $a\in W^s(y)$ by the expansivity of the hyperbolic set.

Because $\alpha(y)=\alpha(a)$ and $\omega(y)=\omega(a)$, we know
that $\alpha(y)\cap \alpha(y)\neq\emptyset$. By the pseudo orbit
tracing property, we can obtain the periodic points $q_n$ with
orbits in $W_1$ such that $q_n\to y$ as $n\to\infty$. It is
obvious that $a$ is close to $\bigcup_{n}W^s(q_n)$. Because
$\{q_n\}$ are pairwise homoclinically related, we know that
$\overline{\bigcup_{n}W^s(Orb(q_n))}=\overline{W^s( Orb(q_n))}$
for any $n$. Further, we know that $a\in \overline{W^s(P)}$.
Similarly, we have $a\in\overline{W^u(Q)}$. This ends the proof of
the lemma.
\hfill$\square$\\

Here, we  complete the proof of Proposition \ref{hetero1}. By clam
and Lemma \ref{diff}, we can take a resist point $a\in
\overline{W^s(P,g)}\cap\overline{W^u(Q,g)}$. Then we can perform a
perturbation in a tube of $a$ such that $W^s(P)\cap W^u(Q)\neq
\emptyset$ and maintain the existing $W^u(P)\cap W^s(Q)\neq
\emptyset$. Thus, we can obtain a heterodimensional cycle.
\hfill$\square$\\

\noindent{\bf End of the proof of Theorem B.} Since $f$ has the
$C^1$ robustly inverse pseudo orbit tracing property on $\Lambda$,
by Lemma \ref{lstar}, $f\in\mathcal{F}(\Lambda).$ We assume that a
transitive set $\Lambda$ is not hyperbolic for $f$. Since
$\Lambda$ is not hyperbolic for $f$, we have a minimally
nonhyperbolic set $K\subset\Lambda$. By Proposition
\ref{nonsimple}, $K$ cannot be a nonsimple-type minimally
nonhyperbolic set. Thus $K$ is a simple-type minimally
nonhyperbolic set of $f$. Then  by Proposition \ref{hetero1},
there is $g$ $C^1$ close to $f$ such that $g$ has a
heterodimensional cycles in $U,$ where $U$ is a locally maximal
neighborhood of $\Lambda.$
  From Proposition \ref{nohetero}, we know that for any $g$ $C^1$
  close to $f$ such that $g$ has no a heterodimensional cycles in
  $U$. Therefore, we can see that $\Lambda$ does not admit the simple-type nonhyperbolic set.
Hence, $\Lambda$ should be a hyperbolic set for $f$.
\hfill$\square$

\bigskip

\bigskip
\noindent {\bf Acknowledgement.} The author would like thanks to
X. Wen for his  valuable comments and suggestions.  This work is
supported by Basic Science Research Program through the National
Research Foundation of Korea (NRF) funded by the Ministry of
Science, ICT \& Future Planning (No. 2017R1A2B4001892 and
2020R1F1A1A01051370).


\end{document}